\documentclass[leqno,11pt]{article}
\usepackage[spanish]{babel}
\usepackage[ansinew]{inputenc}
\usepackage{amscd}
\usepackage{amsthm}
\usepackage{amssymb}
\usepackage{amsmath}
\usepackage{graphics}
\usepackage{graphicx}
\usepackage{verbatim}
\usepackage{color}

\newtheorem{mydef}{Definition}
\newtheorem{myexa}{Example}

\newtheorem{mytheo}{Theorem}

\newtheorem*{myrem}{{\em Remark}}
\newtheorem*{myprem}{{\em Remarks}}
\newtheorem*{myproo}{{\em Proof}}

\newcommand{\pt}{\mbox{$\succ$\hspace{-1ex}$\longrightarrow$}}

\decimalpoint

\begin{document}

\begin{center}
	{\sc On Bayesian Estimation of Densities and Sampling Distributions: the Posterior Predictive Distribution as the Bayes Estimator}\vspace{2ex}\\
	A.G. Nogales\vspace{2ex}\\
	Dpto. de Matem\'aticas, Universidad de Extremadura\\
	Avda. de Elvas, s/n, 06006--Badajoz, SPAIN.\\
	e-mail: nogales@unex.es
\end{center}
\vspace{.4cm}
\begin{quote}
	\hspace{\parindent} {\small {\sc Abstract.} Optimality results for two outstanding Bayesian estimation problems are given in this paper: the estimation of the sampling distribution for the squared total variation function and the estimation of the density for the $L^1$-squared loss function. The posterior predictive distribution provides the solution to these problems. Some examples are presented to illustrate it. The Bayesian estimation problem of a distribution function is also addressed. Consistency of the estimator of the density is proved. 		 
	}
\end{quote}

\vfill
\begin{itemize}
	\item[] \hspace*{-1cm} {\em AMS Subject Class.} (2020): {\em Primary\/} 62F15, 62G07
	{\em Secondary\/} 62C10
	\item[] \hspace*{-1cm} {\em Key words and phrases: } Bayesian density estimation, posterior predictive distribution.
\end{itemize}

\newpage

\section{Introduction and basic definitions}

The posterior predictive distribution has been introduced in the literature to predict the distribution of a future observation from the previous random sample (see Gelman et al. (2014) or Boldstat (2004)), but it has also been used in other areas such as model selection, testing for discordancy, goodness of fit, perturbation analysis or classification. Other fields of application are presented in Geisser (1993) and Rubin (1984). 

In the next pages, the problems of estimation of a density or a probability measure is considered under the Bayesian point of view. These problems are addressed in a number of previous books and papers such as Ghosh et al. (2003, Ch. 5), Lijoi et al. (2010, sect. 3.4), Lo (1984) or Ferguson (1983), to mention just a few. Ghosal et al. (2017), p. 121, contains a brief historical review on Bayesian density estimation. 
Here, the posterior predictive density appears as the optimal estimator of the density for the $L^1$-squared loss function and an arbitrary prior distribution. In fact, the posterior predictive distribution is the optimal estimator of the probability measures $P_\theta$ for the squared total variation loss function. These results come to endorse the applications of the posterior predictive distribution above mentioned.

Several examples are presented in Section 4 to illustrate the results. Gelman et al. (2014) contains many other examples (and exercises) of determination of the posterior predictive distribution. But in practice, the explicit evaluation of the posterior predictive distribution could be cumbersome and its simulation may become preferable.  Gelman et al. (2014) is a good reference also for such simulation methods.

In what follows we will place ourselves in a general framewok for the Bayesian inference, as is described in Barra (1971), for instance, covering simultaneously the discrete, continuous, univariate and multivariate cases. 

First, let us briefly recall some basic concepts about Markov kernels. In the next,  $(\Omega,\mathcal A)$, $(\Omega_1,\mathcal A_1)$ and so on will denote measurable spaces. 

\begin{mydef}\rm  1) (Markov kernel) A Markov kernel
	$M_1:(\Omega,\mathcal A)\pt   (\Omega_1,\mathcal A_1)$ is a map $M_1:\Omega\times\mathcal A_1\rightarrow[0,1]$ such that: 	(i) $\forall \omega\in\Omega$, $M_1(\omega,\cdot)$ is a  probability
	measure on 	$\mathcal A_1$; (ii) $\forall A_1\in\mathcal A_1$, $M_1(\cdot,A_1)$ is $\mathcal A$-measurable.\par
2)	(Image of a Markov kernel) The image (or {\it probability
		distribution}) of a Markov kernel $M_1:(\Omega,\mathcal A,P)\pt
	(\Omega_1,\mathcal A_1)$ on a probability space is the probability
	measure  $P^{M_1}$ on $\mathcal A_1$ defined by
	$P^{M_1}(A_1):=\int_{\Omega}M_1(\omega,A_1)\,dP(\omega)$.
\par
3)  (Composition of Markov kernels) Given two Markov kernels
	$M_1:(\Omega_1,\mathcal A_1)\pt (\Omega_2,\mathcal A_2)$ and
	$M_2:(\Omega_2,\mathcal A_2)\pt (\Omega_3,\mathcal A_3)$, its composition  is defined 	as the Markov kernel $M_2M_1:(\Omega_1,\mathcal A_1)\pt
	(\Omega_3,\mathcal A_3)$ given by
	$$M_2M_1(\omega_1,A_3)=\int_{\Omega_2}M_2(\omega_2,A_3)M_1(\omega_1,d\omega_2).
	$$
\end{mydef}

\begin{myprem}\rm  1) (Markov kernels as extensions of the concept of random variable) The concept of Markov kernel extends the concept of random variable (or measurable map). A random variable $T_1:(\Omega,\mathcal 	A,P)\rightarrow(\Omega_1,\mathcal A_1)$ will be identified with the Markov kernel $M_{T_1}:(\Omega,\mathcal A,P)\pt   (\Omega_1,\mathcal
	A_1)$ defined by $M_{T_1}(\omega,A_1)=\delta_{T_1(\omega)}(A_1)=I_{A_1}(T_1(\omega))$,
	where $\delta_{T_1(\omega)}$ denotes the Dirac measure -the
	degenerate distribution- at the point $T_1(\omega)$, and $I_{A_1}$ is
	the indicator function of the event $A_1$. In particular, the probability distribution $P^{M_{T_1}}$ of $M_{T_1}$ coincides with the probability distribution $P^{T_1}$ of $T_1$ defined as $P^{T_1}(A_1):=P(T_1\in A_1)$\par	
2) Given a Markov kernel $M_1:(\Omega_1,\mathcal A_1)\pt (\Omega_{2},\mathcal A_{2})$ and a random variable $X_2:(\Omega_2,\mathcal A_2)\rightarrow (\Omega_{3},\mathcal A_{3})$, we have that  $M_{X_2}M_1(\omega_1,A_3)=M_1(\omega_1,X_2^{-1}(A_3))=
	M_1(\omega_1,\cdot)^{X_2}(A_3).$ We write $X_2M_1:=M_{X_2}M_1$.
	$\Box$ \end{myprem}

Let $(\Omega,\mathcal A,\{P_\theta\colon \theta\in(\Theta,\mathcal T,Q)\})$ be a Bayesian statistical experiment where $Q$ is the prior distribution, a probability measure on the measurable space $(\Theta,\mathcal T)$. $(\Omega,\mathcal A)$ is the sample space and $(\Theta,\mathcal T)$ is the parameter space. 

When needed, we shall suppose that $P_\theta$ has a density (or Radon-Nikodym derivative) $p_\theta$ with respect to a $\sigma$-finite measure $\mu$ on $\mathcal A$ and that the likelihood function $\mathcal L:(\omega,\theta)\in(\Omega\times\Theta,\mathcal A\otimes\mathcal T)\rightarrow \mathcal L(\omega,\theta):=p_\theta(\omega)$ is measurable. So we have a Markov kernel $P:(\Theta,\mathcal T)\pt (\Omega,\mathcal A)$ defined by $P(\theta,A):=P_\theta(A)$. Let $P^*:(\Omega,\mathcal A)\pt (\Theta,\mathcal T)$ the Markov kernel determined by the posterior distributions. In fact, if we denote by $\Pi$ the only probability measure  on $\mathcal A\otimes \mathcal T$ such that 
\begin{gather*}
	\Pi(A\times T)=\int_TP_\theta(A)dQ(\theta), \quad A\in\mathcal A, T\in\mathcal T,\tag{1}
\end{gather*}
then $P^*$ is defined in such a way that 
\begin{gather*}\Pi(A\times T)=\int_AP^*_\omega(T)d\beta_Q^*(\omega),\quad A\in\mathcal A, T\in\mathcal T,\tag{2}
\end{gather*}
where $\beta_Q^*$ denotes the so called prior predictive probability, defined by 
$$\beta_Q^*(A)=\int_\Theta P_\theta(A)dQ(\theta), \quad A\in\mathcal A.
$$
In other terms, $\beta_Q^*=Q^P$, the probability distribution of the Markov kernel $P$ with respect to the prior distribution $Q$. 

The probability measure $\Pi$ integrates all the basic ingredients of the Bayesian model, and these ingredients can be essentially derived from $\Pi$, something that would allow us to identify the Bayesian model as the probability space $(\Omega\times\Theta,\mathcal A\otimes\mathcal T,\Pi)$ (so is done, for instance, in Florens et al. (1990)). 

It is well known that, for $\omega\in\Omega$, the posterior density  with respect to the prior distribution is proportional to the likelihood. Namely
$$p^*_\omega(\theta):=\frac{dP^*_\omega}{dQ}(\theta)=C(\omega)p_\theta(\omega),
$$
where $C(\omega)=\left[\int_\Theta p_\theta(\omega)dQ(\theta)\right]^{-1}$. 

This way we obtain a statistical experiment $(\Theta,\mathcal T,\{P_\omega^*\colon \omega\in\Omega\})$ on the parameter space $(\Theta,\mathcal T)$. 
 We can reconsider the Markov kernel $P$ defined on this statistical experiment
$$P:(\Theta,\mathcal T,\{P_\omega^*\colon \omega\in\Omega\})\pt (\Omega,\mathcal A).
$$
Since $\big(P_\omega^*\big)^P(A)=\int_\Theta P_\theta(A)dP_\omega^*(\theta)$, for $A\in\mathcal A$, it is  called the posterior predictive distribution on $\mathcal A$ given $\omega$, and the statistical experiment image of $P$ is
$$\big(\Omega,\mathcal A,\big\{\big(P_\omega^*\big)^P\colon \omega\in\Omega\big\}\big).
$$
Note that, given $\omega\in\Omega$, according to Fubini's Theorem,
\begin{gather*}\begin{split}
		\big(P_\omega^*\big)^P(A)&=\int_\Theta P_\theta(A)dP_\omega^*(\theta)=\int_\Theta \int_A p_\theta(\omega')d\mu(\omega')p_\omega^*(\theta)dQ(\theta)\\
		&=\int_A\int_\Theta p_\theta(\omega')p_\omega^*(\theta)dQ(\theta)d\mu(\omega').
\end{split}\end{gather*}

So, the posterior predictive density is
$$\frac{d\big(P_\omega^*\big)^P}{d\mu}(\omega')=\int_\Theta p_\theta(\omega')p_\omega^*(\theta)dQ(\theta).
$$

If we consider the composition of the Markov kernels $P^*$ and $P$:
$$(\Omega,\mathcal A)\overset{P^*}{\pt}(\Theta,\mathcal T)\overset{P}{\pt}(\Omega,\mathcal A),
$$
defined by 
\begin{gather*}
	PP^*(\omega,A):=\int_\Theta P_\theta(A)dP^*_\omega(\theta)=
\int_A\int_\Theta p_\theta(\omega')p^*_\omega(\theta)dQ(\theta)d\mu(\omega'),\tag{3}
\end{gather*}
we have that
$$\frac{dPP^*(\omega,\cdot)}{d\mu}(\omega')=\int_\Theta p_\theta(\omega')p^*_\omega(\theta)dQ(\theta).
$$

Notice that $PP^*(\omega,\cdot)=\big(P^*_\omega\big)^P$.

\begin{myprem}\rm 1) Analogously, if we consider the Markov kernel 
	$$P^*:(\Omega,\mathcal A,\{P_\theta\colon \theta\in\Theta\})\pt (\Theta,\mathcal T,\{\big(P_\theta\big)^{P^*}\colon\theta\in\Theta\}),$$ 
	we have that 
	$$\big(P_\theta\big)^{P^*}(T)=\int_\Omega P^*_\omega(T)dP_\theta(\omega)=\int_T\int_\Omega p_\omega^*(\theta')p_\theta(\omega)d\mu(\omega)dQ(\theta')
	$$
	and
	$$\frac{d\big(P_\theta\big)^{P^*}}{dQ}(\theta')=
	\int_\Omega p_\omega^*(\theta')p_\theta(\omega)d\mu(\omega).
	$$
	Notice that $\big(\beta_Q^*\big)^{P^*}(T)=\int_\Omega P^*_\omega(T)d\beta_Q^*(\omega)=Q(T)$, i.e.,
	$\big(\beta_Q^*\big)^{P^*}=Q$.
	
	If we consider the composition of the Markov kernels $P$ and $P^*$:
	$$(\Theta,\mathcal T)\overset{P}{\pt}(\Omega,\mathcal A)\overset{P^*}{\pt}(\Theta,\mathcal T),
	$$
	defined by
	$$P^*P(\theta,T)=\int_\Omega P^*_\omega(T)dP_\theta(\omega)=
	\int_T\int_\Omega p^*_\omega(\theta')p_\theta(\omega)d\mu(\omega)dQ(\theta'),
	$$
	we obtain
	$$\frac{dP^*P(\theta,\cdot)}{dQ}(\theta')=\int_\Omega p^*_\omega(\theta')p_\theta(\omega)d\mu(\omega).
	$$
	
	Notice that $P^*P(\theta,\cdot)=P_\theta^{P^*}$. As a consequence of Fubini's Theorem, we get
	$$Q^{P^*P}=Q.$$\par
	2) Because of (1), we introduce the notation $\Pi:=P\otimes Q$. So, (2) reads as $\Pi:=\beta_Q^*\otimes P^*$. Hence, after observing $\omega\in\Omega$, replacing the prior distribution $Q$ by the posterior distribution $P_\omega^*$, we get the probability distribution $\Pi_\omega:=P\otimes P_\omega^*$ on $\mathcal A\otimes\mathcal T$. According to (3),  $PP^*(\omega,A)=\Pi_\omega(A\times\Theta)=\Pi_\omega^I(A)$ where $I(\omega,\theta)=\omega$. This way the posterior predictive distribution $\big(P^*_\omega\big)^P$ given $\omega$ appears as the marginal $\Pi_\omega$-distribution on $\Omega$. $\quad\Box$	
\end{myprem}

\section{Bayesian estimation of probabilities, sampling distributions and densities}

 According to Bayesian philosophy, given $A\in\mathcal A$, a natural estimator of $f_A(\theta):=P_\theta(A)$ is the posterior mean of $f_A$, which coincides with  the posterior predictive probability of $A$, $T(\omega):=\big(P_\omega^*\big)^P(A)$. In fact, this is the Bayes estimator of $f_A$ (see Theorem 1.(i)). 
	
So, the posterior predictive distribution $\big(P_\omega^*\big)^P$ appears as the natural Bayesian estimator of the probability distribution $P_\theta$. 
	
To estimate probability measures, the squared total variation loss function
	$$W_1(Q,P):=\sup_{A\in\mathcal A}|Q(A)-P(A)|^2,
	$$
	will be considered.
	An estimator of $f(\theta):=P_\theta$ is a Markov kernel $M:(\Omega,\mathcal A)\pt (\Omega,\mathcal A)$ so that, being observed $\omega\in\Omega$, $M(\omega,\cdot)$ is a probability measure on $\mathcal A$ which is considered as an estimation of $f$. We wonder if the Bayes mean risk of the estimator  $M^*:=\big(P^*\big)^P$ is less than that of any other estimator $M$ of $f$, i.e., we wonder if
	$$\int_{\Omega\times\Theta}\sup_{A\in\mathcal A}|\big(P_\omega^*\big)^P(A)-P_\theta(A)|^2d\Pi(\omega,\theta)\le
	\int_{\Omega\times\Theta}\sup_{A\in\mathcal A}|M(\omega,A)-P_\theta(A)|^2d\Pi(\omega,\theta).$$
	Theorem 1.(ii) below gives the answer.

An estimator of the density $p_\theta$ on $(\Omega,\mathcal A,\{P_\theta\colon \theta\in(\Theta,\mathcal T,Q)\})$ is a measurable map $m:(\Omega^2,\mathcal A^2)\longrightarrow \mathbb R$ in such a way that, being observed $\omega\in\Omega$, the map $\omega'\mapsto m(\omega,\omega')$ is an estimation of $p_\theta$. 

It is well known (see Ghosal et al. (2017), p. 126) that, given two probability measures $Q$ and $P$ on $(\Omega,\mathcal A)$ having densities $q$ and $p$ with respect to a $\sigma$-finite measure $\mu$, 
$$\sup_{A\in\mathcal A}|Q(A)-P(A)|=\frac12\int|q-p|d\mu.
$$

So the Bayesian estimation of the sampling distribution $P_\theta$ for the squared total variation loss function corresponds to the Bayesian estimation of its density $p_\theta$ for the $L^1$-squared loss function
$$W'_1(q,p):=\big(\textstyle\int|q-p|\,d\mu\big)^2,
$$

The next Theorem also solves the estimation problem of the density.

\begin{mytheo}\rm Let $(\Omega,\mathcal A,\{P_\theta\colon \theta\in(\Theta,\mathcal T,Q)\})$ be a Bayesian statistical experiment dominated by a $\sigma$-finite measure $\mu$, where the $\sigma$-field $\mathcal A$ is supposed to be separable.   We suppose that the likelihood function $\mathcal L(\omega,\theta):=p_\theta(\omega)=dP_\theta(\omega)/d\mu$ is $\mathcal A\otimes\mathcal T$-measurable.

(i) Given $A\in\mathcal A$, the posterior predictive probability $\big(P_\omega^*\big)^P(A)$ of $A$ is the Bayes estimator of the probability $P_\theta(A)$ of $A$ for the squared error loss function $$W(x,\theta):=(x-P_\theta(A))^2.$$
Moreover, if $X$ is a real statistics with finite mean, its posterior predictive mean
$$E_{(P_\omega^*)^P}(X)=\int_\Theta\int_\Omega X(\omega')dP_\theta(\omega')dP_\omega^*(\theta)$$
is the Bayes estimator  of $E_\theta(X)$.

(ii) The posterior predictive distribution $\big(P_\omega^*\big)^P$ is the Bayes estimator of the sampling distribution $P_\theta$ for the squared total variation loss function $$W_1(P,Q):=\sup_{A\in\mathcal A}|P(A)-Q(A)|^2.
$$

(iii) The posterior predictive density 
$$b^*_{Q,\omega}(\omega'):=\frac{d\big(P_\omega^*\big)^P}{d\mu}(\omega')=\int_\Theta p_\theta(\omega')p_\omega^*(\theta)dQ(\theta).
$$
is the Bayes estimator of the density $p_\theta$ for the $L^1$-squared loss function
$$W'_1(p,q):=\left(\int_\Omega |p-q|d\mu\right)^2.
$$
	\end{mytheo}

\section{Bayesian estimation of sampling distributions and densities from a sample}

	More generally, an estimator of $f(\theta):=P_\theta$ from a sample of size $n$ of this distribution is a Markov kernel 
	$$M_n:(\Omega^n,\mathcal A^n)\pt (\Omega,\mathcal A).
	$$
	Let us consider the Markov kernel
	$$P^n:(\Theta,\mathcal T)\pt (\Omega^n,\mathcal A^n)
	$$
	defined by $P^n(\theta,A)=P_\theta^n(A)$, $A\in\mathcal A^n$, $\theta\in\Theta$. We write
	$\Pi_n:=P^n\otimes Q,$
	so that 
	$$\Pi_n(A\times T)=\int_TP_\theta^n(A)dQ(\theta),\quad A\in\mathcal A^n,T\in\mathcal T.
	$$
	The corresponding prior predictive distribution is
	$$\beta_{Q,n}^*(A)=\int_\Theta P^n_\theta(A)dQ(\theta)=\Pi_n^I(A),
	$$
	where $I(\omega,\theta)=\omega$ for $\omega\in\Omega^n$. Let us write $I_i(\omega)=\omega_i$ and $\hat I_i(\omega,\theta)=\omega_i$, for $\omega\in\Omega^n$ and $i=1,\dots,n$. Hence
	$$\big(\beta_{Q,n}^*\big)^{I_i}(A_i)=\int_\Theta P_\theta(A_i)dQ(\theta)=\beta_Q^*(A_i),
	$$
	and
	$$\Pi_n^{\hat I_i}(A_i\times T)=\int_TP_\theta(A_i)dQ(\theta),
	$$
	so 
	$$\big(\beta_{Q,n}^*\big)^{I_i}=\beta_Q^*,\qquad\mbox{and}\qquad \Pi_n^{\hat I_i}=\Pi.
	$$
	Denoting $J(\omega,\theta)=\theta$, the posterior distribution $P^*_{\omega,n}:=\Pi_n^{J|I=\omega}$, $\omega\in\Omega^n$, is defined in such a way that
	$$\Pi_n(A\times T)=\int_A P^*_{\omega,n}(T)d\beta_{Q,n}^*(\omega).
	$$
	The $\mu^n$-density of $P_\theta^n$ is	$$p_{\theta,n}(\omega):=\frac{dP_\theta^n}{d\mu^n}(\omega)=\prod_{i=1}^np_\theta(\omega_i)\quad \mbox{for\ } \omega=(\omega_1,\dots,\omega_n)\in\Omega^n.
	$$
	The posterior density given $\omega\in\Omega^n$ is of the form
	$$p^*_{\omega,n}(\theta):=\frac{dP^*_{\omega,n}}{dQ}(\theta)\propto p_{\theta,n}(\omega).
	$$
	According to Theorem 1.(ii), the Markov kernel 
	$$\big(P^*_{n}\big)^{P^n}:(\Omega^n,\mathcal A^n)\pt (\Omega^n,\mathcal A^n)
	$$
	defined by
	$$\big(P^*_{n}\big)^{P^n}(\omega,A):=\big(P^*_{\omega,n}\big)^{P^n}(A)=\int_\Theta P_\theta^n(A)dP^*_{\omega,n}(\theta),
	$$
	is the Bayes estimator of the product probability measure $f_n(\theta):=P_\theta^n$. That is to say
	$$\int_{\Omega\times\Theta}\sup_{A\in\mathcal A^n}|\big(P_{\omega,n}^*\big)^{P^n}(A)-P_\theta^n(A)|^2d\Pi_n(\omega,\theta)\le
	\int_{\Omega\times\Theta}\sup_{A\in\mathcal A^n}|M(\omega,A)-P_\theta^n(A)|^2d\Pi_n(\omega,\theta),$$
	for every estimator $M:(\Omega^n,\mathcal A^n)\pt (\Omega^n,\mathcal A^n)$ of $P_\theta^n$.

The next theorem shows how marginalizing the posterior predictive distribution $\big(P_{\omega,n}^*\big)^{P^n}$ we can get the Bayes estimator of the sampling probability measure $P_\theta$ or its density.

\begin{mytheo} \rm (Bayesian density estimation from a sample of size $n$) Let $(\Omega,\mathcal A,\{P_\theta\colon \theta\in(\Theta,\mathcal T,Q)\})$ be a Bayesian statistical experiment dominated by a $\sigma$-finite measure $\mu$, where the $\sigma$-field $\mathcal A$ is supposed to be separable.   We suppose that the likelihood function $\mathcal L(\omega,\theta):=p_\theta(\omega)=dP_\theta(\omega)/d\mu$ is $\mathcal A\otimes\mathcal T$-measurable. Let $n\in\mathbb N$. All the estimation problems below are referred to the product Bayesian statistical experiment $(\Omega^n,\mathcal A^n,\{P_\theta^n\colon \theta\in(\Theta,\mathcal T,Q)\})$ corresponding to a $n$-sized sample of the observed unknown distribution. Let $I_1(\omega_1,\dots,\omega_n):=\omega_1$. 
	
	(i) Given $A\in\mathcal A$,  $$\left[\big(P_{\omega,n}^*\big)^{P^n}\right]^{I_1}(A)$$  is the Bayes estimator of the probability $P_\theta(A)$ of $A$ for the squared error loss function $$W(x,\theta):=(x-P_\theta(A))^2.$$

	(ii) The distribution $$\left[\big(P_{\omega,n}^*\big)^{P^n}\right]^{I_1}$$ 
	of the projection $I_1$ under the posterior predictive probability $\big(P_{\omega,n}^*\big)^{P^n}$ is the Bayes estimator of the sampling distribution $P_\theta$ for the squared total variation loss function $$W_1(P,Q):=\sup_{A\in\mathcal A}|P(A)-Q(A)|^2.
	$$
	
	(iii) The marginal posterior predictive density  
	$$b^*_{Q,\omega,n}(\omega'):=\frac{d\left[\big(P_{\omega,n}^*\big)^{P^n}\right]^{I_1}}{d\mu}(\omega')=\int_\Theta p_\theta(\omega')p_{\omega,n}^*(\theta)dQ(\theta).
	$$
	is the Bayes estimator of the density $p_\theta$ for the $L^1$-squared loss function
	$$W'_1(p,q):=\left(\int_\Omega |p-q|d\mu\right)^2.
	$$
\end{mytheo}

We end this section with a remark that address the problem of estimating a real distribution function. 

\begin{myrem}\rm (Bayesian estimation of a distribution function) When $P_\theta$ is a probability distribution on the line, we may be interested in the estimation of its distribution function $F_\theta(t):=P_\theta(]-\infty,t])$. An estimator of such a distribution function is a map 
	$$F:(x,t)\in\mathbb R^n\times\mathbb R\longmapsto F(x,t):=M(x,]-\infty,t])
	$$
	for a Markov kernel $M:(\mathbb R^n,\mathcal R^n)\pt (\mathbb R,\mathcal R)$, where $\mathcal R$ denotes the Borel $\sigma$-field on $\mathbb R$.
	
	Accordig to the previous results, given $t\in\mathbb R$,
	$$F_{x}^*(t):=\left[\big(P_{x,n}^*\big)^{P^n}\right]^{I_1}(]-\infty,t])=\int_{-\infty}^t\int_\Theta p_{\theta,n}(y)\cdot p^*_{x,n}(\theta)dQ(\theta)d\mu^n({y})
	$$
	is the Bayes estimator of $F_\theta(t)$ for the squared error loss function. So
	$$\int_{\mathbb R^n\times\Theta}|F_{x}^*(t)-F_\theta(t)|^2d\Pi(x,\theta)\le
	\int_{\mathbb R^n\times\Theta}|F(x,t)-F_\theta(t)|^2d\Pi(x,\theta)
	$$
	for any other estimator $F$ of $F_\theta$.  Since 
	$$\sup_{t\in\mathbb R}|F(x,t)-F_\theta(t)|=\sup_{r\in\mathbb Q}|F(x,r)-F_\theta(r)|
	$$
	we have that, given $(x,\theta)\in\mathbb R^n\times\Theta$ and $k\in\mathbb N$, there exists $r_k\in\mathbb Q$ such that 
	$$C(x,\theta)-\frac1k\le|F_{x}^*(r_k)-F_\theta(r_k)|,
	$$
	where $C(x,\theta):=\sup_{t\in\mathbb R}|F_{x}^*(t)-F_\theta(t)|^2$, and hence (see the remark at the end of Section \nolinebreak 5)
	\begin{gather*}\begin{split}
			\int_{\mathbb R^n\times\Theta}C(x,\theta)d\Pi(x,\theta)&\le \int_{\mathbb R^n\times\Theta}|F_{x}^*(r_k)-F_\theta(r_k)|^2d\Pi(x,\theta)+\frac1k\\
			&\le \int_{\mathbb R^n\times\Theta}\sup_{t\in\mathbb R}|F(x,t)-F_\theta(t)|^2d\Pi(x,\theta)+\frac1k.
		\end{split}
	\end{gather*}
	We have proved that the posterior predictive distribution function $F_{x}^*$ is the Bayes estimator of the distribution function $F_\theta$ for the $L^{\infty}$-squared loss function
	$$W''(F,G)=\big(\sup_{t\in\mathbb R}|F(t)-G(t)|\big)^2.\quad\Box
	$$
	
\end{myrem}

\section{Consistency of the Bayesian estimator of the density}

First we adapt the framework to an asymptotic context. 
Let
$$(\Omega,\mathcal A,\{P_\theta\colon \theta\in(\Theta,\mathcal T,Q)\})
$$ 
be a Bayesian experiment and consider the product Bayesian experiment
$$(\Omega^{\mathbb N},\mathcal A^{\mathbb N},\{P_\theta^{\mathbb N}\colon \theta\in(\Theta,\mathcal T,Q)\})
$$
corresponding to an infinite sample of the unknown distribution $P_\theta$. Let us write
$$I(\omega,\theta):=\omega,\quad J(\omega,\theta):=\theta,\quad I_n(\omega,\theta):=\omega_n\quad \mbox{and}\quad I_{(n)}(\omega):=\omega_{(n)}:=(\omega_1,\dots,\omega_n),
$$
where $n\in\mathbb N$. 

We suppose that  $P^{\mathbb N}(\theta,A):=P_\theta^{\mathbb N}(A)$ is a  Markov kernel. 
Let
$$\Pi_{\mathbb N}:=P^{\mathbb N}\otimes Q
$$
i.e., 
$$\Pi_{\mathbb N}(A\times T)=\int_TP_\theta^{\mathbb N}(A)dQ(\theta),\quad A\in\mathcal A,\; T\in\mathcal T.
$$

Being $Q:=\Pi_{\mathbb N}^{J}$,  $P_\theta^{\mathbb N}$ is a version of the conditional distribution (regular conditional probability) $\Pi_{\mathbb N}^{I|J=\theta}$. Analogously, $P_\theta^n$ is a version of the conditional distribution $\Pi_{\mathbb N}^{I_{(n)}|J=\theta}$.

Let $\beta_{Q,\mathbb N}^*:=\Pi_{\mathbb N}^{I}$, the prior predictive distribution in $\Omega^{\mathbb N}$. Similarly we write  $\beta_{Q,n}^*:=\Pi_{\mathbb N}^{I_{(n)}}$ for the prior predictive distribution in $\Omega^n$. So, the posterior distribution is $P^*_{\omega,\mathbb N}:=\Pi_{\mathbb N}^{J|I=\omega}$, in such a way that 
$$\Pi_{\mathbb N}(A\times T)=\int_TP_\theta^{\mathbb N}(A)dQ(\theta)=\int_AP^*_{\omega,\mathbb N}(T)d\beta_{Q,\mathbb N}^*(\omega),\quad A\in\mathcal A,\; T\in\mathcal T.
$$
Denote
$P^*_{\omega_{(n)},n}:=\Pi_{\mathbb N}^{J|I_{(n)}=\omega_{(n)}}$ for $\omega_{(n)}\in\Omega^n$.

Write ${P^*_{\omega_{(n)},n}}^{\!\!\!\!\!\!\!\!\!\!P\,\,\,\,\,}$ for the posterior predictive distribution given $\omega_{(n)}\in\Omega^n$, defined for $A\in\mathcal A$ as
$${P^*_{\omega_{(n)},n}}^{\!\!\!\!\!\!\!\!\!\!P\,\,\,\,\,}(A)=
\int_\Theta P_\theta(A)dP^*_{\omega_{(n)},n}(\theta).
$$

In the dominated case we can assume without loss of generality that the dominating measure $\mu$ is a probability measure (because of (1) below). We write $p_\theta=dP_\theta/d\mu$. 

We have that, far all $n$ and every event $A\in\mathcal A$, 
\begin{gather*}\begin{split}
		{P^*_{\omega_{(n)},n}}^{\!\!\!\!\!\!\!\!\!\!P\,\,\,\,\,}(A)&=
		\int_\Theta P_\theta(A)dP^*_{\omega_{(n)},n}(\theta)=
		\int_\Theta \int_A p_\theta(\omega')d\mu(\omega')dP^*_{\omega_{(n)},n}(\theta)
		\\&
		=
		\int_A\int_\Theta  p_\theta(\omega')dP^*_{\omega_{(n)},n}(\theta)d\mu(\omega')
\end{split}\end{gather*}
which proves that
$${p^*_{\omega_{(n)},n}}^{\!\!\!\!\!\!\!\!\!\!P\,\,\,\,\,}(\omega'):=\int_{\Theta} p_\theta(\omega')dP^*_{\omega_{(n)},n}(\theta)
$$
is a $\mu$-density of  ${P^*_{\omega_{(n)},n}}^{\!\!\!\!\!\!\!\!\!\!P\,\,\,\,\,}$.

In the same way
$${p^*_{\omega,\mathbb N}}^{\!\!\!\!P}(\omega'):=\int_{\Theta} p_\theta(\omega')dP^*_{\omega,\mathbb N}(\theta)
$$
is a $\mu$-density of ${P^*_{\omega,\mathbb N}}^{\!\!\!\!P}$, the posterior predictive distribution given $\omega\in\Omega^{\mathbb N}$. 

It has been previously proved by the author (see Nogales (2020)) that the posterior predictive distribution
${P^*_{\omega_{(n)},n}}^{\!\!\!\!\!\!\!\!\!\!P\,\,\,\,\,}$
(resp. the posterior predictive density ${p^*_{\omega_{(n)},n}}^{\!\!\!\!\!\!\!\!\!\!P\,\,\,\,\,}$) is the Bayes estimator of the sampling distribution $P_\theta$ (resp. the density $p_\theta$) for the squared variation total (resp. the squared $\mathbb L^1$) loss function  in the product experiment $(\Omega^n,\mathcal A^n,\{P_\theta^n\colon \theta\in(\Theta,\mathcal T,Q)\})$. Analogously, the posterior predictive distribution 
${P^*_{\omega,\mathbb N}}^{\!\!\!\!P}$ (resp. the posterior predictive density ${p^*_{\omega,\mathbb N}}^{\!\!\!\!P}$) is the Bayes estimator of the sampling distribution $P_\theta$ (resp. the density $p_\theta$) for the squared variation total (resp. the squared $\mathbb L^1$) loss function  in the product experiment $(\Omega^{\mathbb N},\mathcal A^{\mathbb N},\{P_\theta^{\mathbb N}\colon \theta\in(\Theta,\mathcal T,Q)\})$.

It is well known that 
$$\sup_{A\in\mathcal A}\left|{P^*_{\omega_{(n)},n}}^{\!\!\!\!\!\!\!\!\!\!P\,\,\,\,\,}(A)-P_\theta(A)\right|=\frac12\int_{\Omega}\left|{p^*_{\omega_{(n)},n}}^{\!\!\!\!\!\!\!\!\!\!P\,\,\,\,\,}-p_\theta\right|d\mu.\qquad (1)
$$

We wonder if the Bayes risk of the Bayes estimator ${P^*_{\omega_{(n)},n}}^{\!\!\!\!\!\!\!\!\!\!P\,\,\,\,\,}$ of the sampling distribution $P_\theta$ goes to zero when  $n\to\infty$, i.e., if 
$$\lim_n\int_{\Omega^{\mathbb N}\times\Theta}\sup_{A\in\mathcal A}\left|{P^*_{\omega_{(n)},n}}^{\!\!\!\!\!\!\!\!P\,\,\,\,}(A)-P_\theta(A)\right| d\Pi_{\mathbb N}(\omega,\theta)=0.
$$

In terms of densities, the question is whether  the Bayes risk of the Bayes estimator  ${p^*_{\omega_{(n)},n}}^{\!\!\!\!\!\!\!\!\!\!P\,\,\,\,\,}$ of the density $p_\theta$ goes to zero when  $n\to\infty$, i.e., if 
$$\lim_n \int_{\Omega^{\mathbb N}\times\Theta}\left(\int_{\Omega}\left|{p^*_{\omega_{(n)},n}}^{\!\!\!\!\!\!\!\!\!\!P\,\,\,\,\,}(\omega')-p_\theta(\omega')\right|d\mu(\omega')\right)^2d\Pi_{\mathbb N}(\omega,\theta)=0.
$$

The following result, consequence of a Theorem of Doob that the reader can find in Ghosal et al. (2017), provides the answers and shows the consistency of the estimator. 

\begin{mytheo}\rm Let $(\Omega,\mathcal A,\{P_\theta\colon \theta\in(\Theta,\mathcal T,Q)\})$ be a Bayesian experiment dominated by a $\sigma$-finite measure $\mu$. Let us suppose that $(\Omega,\mathcal A)$ is a  Borel standar space, that $\Theta$ is a Borel subset of a polish space and $\mathcal T$ is its Borel  $\sigma$-field. Suppose also that the likelihood function  $\mathcal L(\omega,\theta):=p_\theta(\omega)=\frac{dP_\theta}{d\mu}(\omega)$ is $\mathcal A\times\mathcal T$-measurable and the family $\{P_\theta\colon \theta\in\Theta\}$ is identifiable. Then:
	\begin{itemize}
		\item[(a)] The posterior predictive density  ${p^*_{\omega_{(n)},n}}^{\!\!\!\!\!\!\!\!P\,\,\,\,}$ 
		is the Bayes estimator of the density  $p_\theta$ in the product experiment  $(\Omega^n,\mathcal A^n,\{P_\theta^n\colon \theta\in(\Theta,\mathcal T,Q)\})$ for the  squared $\mathbb L^1$ loss function. Moreover the risk function converges to 0 both for the $\mathbb L^1$ loss function and the squared  $\mathbb L^1$ loss function. 
		
		\item[(b)] The posterior predictive distribution  ${P^*_{\omega_{(n)},n}}^{\!\!\!\!\!\!\!\!P\,\,\,\,}$ 
		is the Bayes estimator of the sampling distribution  $P_\theta$ in the product experiment  $(\Omega^n,\mathcal A^n,\{P_\theta^n\colon \theta\in(\Theta,\mathcal T,Q)\})$ for the squared variation total loss function. Moreover the risk function converges to 0 both for the variation total loss function and the squared variation total loss function. 
		
		\item[(c)] The posterior predictive density is a consistent estimator of the density $p_\theta$, i.e.,
		$$\lim_n {p^*_{\omega_{(n)},n}}^{\!\!\!\!\!\!\!\!\!\!P\,\,\,\,\,}(\omega')=p_\theta(\omega'),\quad \mu\times P_\theta^{\mathbb N}-\hbox{a.e.}
		$$
		for $Q$-almost every $\theta\in\Theta$. 
		
	\end{itemize}
\end{mytheo}

\section{Examples}

\begin{myexa}\rm Let $G(\alpha,\beta)$ be the distribution gamma with parameters $\alpha,\beta>0$ and $P_\theta:=G(1,\theta^{-1})$, whose density is $p_\theta(x)=\theta\exp\{-\theta x\}$ for $x>0$. 
	
	So $P_\theta^n$ is the joint distribution of a sample of size $n$ of an exponential distribution of parameter $1/\theta$ and its density is $p_{\theta,n}(x)=\theta^n\exp\{-\theta\sum_ix_i\}$ for $x=(x_1,\dots,x_n)\in\mathbb R_+^n$.  
	
	Consider the prior distribution $Q:=G(1,\lambda^{-1})$ for some known $\lambda>0$. 
	
	Since, for $a>0$,
	$$\int_0^\infty \theta^n\exp\{-a\theta\}d\theta=\frac{n!}{a^{n+1}},
	$$
	we have that the posterior density given $x\in\mathbb R_+^n$ is
	$$p_{x,n}^*(\theta)=\frac{\lambda}{n!}\big(\textstyle\sum_ix_i\big)^{n+1}\theta^n\exp\{-\theta(\lambda+\sum_ix_i)\}.
	$$
	So, denoting by $\mu_n$ the Lebesgue measure on $\mathbb R_+^n$, the density of the posterior predictive probability given $x$ is
	$$\frac{d(P_{x,n}^*)^{P^n}}{d\mu_n}(x')=\int_\Theta p_{\theta,n}(x')\cdot p_{x,n}^*(\theta)\,d\theta=\lambda\frac{(2n)!}{n!}\frac{\big(\textstyle\sum_ix_i\big)^{n+1}}{\big(\lambda+\textstyle\sum_ix'_i+\sum_ix_i\big)^{2n+1}}.
	$$	
	According to the previous results, this is the Bayes estimator of the joint density $p_{\theta,n}$ for the loss function
	$$W'_n(q,p):=\left(\int_{\mathbb R^n}|q-p|d\mu_n\right)^2,
	$$ 
	while the posterior predictive distribution $\big(P_{x,n}^*\big)^{P^n}$ is the Bayes estimator of the sampling distribution $P_\theta^n$ for the squared total variation loss function on $(\Omega^n,\mathcal A^n)$. 
	
	Moreover, the image $M_n^*(x,\cdot):=\left[\big(P_{x,n}^*\big)^{P^n}\right]^{I_1}=I_1\big(P_{x,n}^*\big)^{P^n}$ is the Bayes estimator of the probability distribution $P_\theta$ for the squared total variation on $(\Omega,\mathcal A)$ and its density
	$$x'>0\,\,\longmapsto\,\,\frac{dM_n^*(x,\cdot)}{d\mu_1}(x')=\int_0^\infty p_\theta(x')\cdot p^*_{x,n}(\theta)\,d\theta=\frac{n\lambda\big(\sum_{i=1}^nx_i\big)^{n+1}}{(\lambda+x'+\sum_{i=1}^nx_i)^{n+2}}
	$$
	is the Bayes estimator of the density $p_\theta$ for the $L^1$-squared loss function $W'_1$.$\quad\Box$
\end{myexa}

\begin{myexa}\rm Let $P_\theta$ the normal distribution $N(\theta,\sigma_0^2)$ with unknown mean $\theta\in\mathbb R$ and known variance $\sigma_0^2$. Let $Q:=N(\mu,\tau^2)$ be the prior distribution where the mean $\mu$ and variance $\tau^2$ are known constants. It is well known that the posterior distribution is $P^*_{x,n}=N(m_n(x),s_n^2)$ where
	$$m_n(x)=\frac{n\tau^2\bar x+\sigma_0^2\mu}{n\tau^2+\sigma_0^2}\quad\mbox{and}\quad s_n^2=\frac{\tau^2\sigma_0^2}{n\tau^2+\sigma_0^2}.
	$$
It can be shown  that the distributiion of $I_1$ with respect to the posterior predictive distribution is
$$\left[\big(P^*_{x,n}\big)^{P^n}\right]=N(m_n(x),\sigma_0^2+s_n^2).
$$
For the details, the reader is addressed to Boldstat (2004, p. 185), where the distributiion of $I_1$ with respect to the posterior predictive distribution is referred to as the predcitive distribution for the next observation given the observation $x$. 

So $M_n^*(x,\cdot):=N(m_n(x),\sigma_0^2+s_n^2)$ is the Bayes estimator of the sampling distribution $N(\theta,\sigma_0^2)$ for the squared total variation loss function and the density of $N(m_n(x),\sigma_0^2+s_n^2)$ is the Bayes estimator of the density of $N(\theta,\sigma_0^2)$ for the $L^1$-squared loss function. 
	$\quad\Box$

\end{myexa}
\begin{myexa}\rm Let $P_\theta$ be the Poisson distribution with parameter $\theta>0$ whose probablity function (or density with respect to the counter measure $\mu_1$ on $\mathbb N_0$) is $p_\theta(k)=\exp\{-\theta\}\frac{\theta^k}{k!}$ for $k\in\mathbb N_0$. 
	
	So $P_\theta^n$ is the joint distribution of a sample of size $n$ of a Poisson distribution of parameter $\theta$ and its probability function (or density with respecto to the counter measure $\mu_n$ on $\mathbb N_0^n$) is $p_{\theta,n}(k)=\exp\{-n\theta\}\frac{\theta^{\|k\|_1}}{\prod_{i=1}^n(k_i!)}$ for $k=(k_1,\dots,k_n)\in\mathbb N_0^n$, where $\|k\|_1:=\sum_{i=1}^nk_i$.  
	
	Consider the prior distribution $Q:=G(1,\lambda^{-1})$ for some known $\lambda>0$. 
	
	It is readily shown that the posterior distribution given $k\in\mathbb N_0^n$ is the gamma distribution $G\big(\|k\|_1,\frac1{\lambda+n}\big)$ whose density is
	$$p_{k,n}^*(\theta)=\frac{(\lambda+n)^{\|k\|_1}}{(\|k\|_1)!}\cdot\theta^{\|k\|_1}\exp\{-\theta(\lambda+n)\}.
	$$
	So the probability function of the posterior predictive probability given $k\in\mathbb N_0^n$ is
	$$\frac{d(P_{k,n}^*)^{P^n}}{d\mu_n}(k')=\int_\Theta p_{\theta,n}(k')\cdot p_{k,n}^*(\theta)\,d\theta=\frac{(\|k'\|_1+\|k\|_1)!}{\prod_{i=1}^n(k_i!)\cdot(\|k\|_1)!}\cdot\frac{(\lambda+n)^{\|k\|_1+1}}{(\lambda+2n)^{\|k'\|_1+\|k\|_1+1}}.
	$$	
	According to the previous results, this is the Bayes estimator of the joint density $p_{\theta,n}$ for the loss function
	$$W'_n(q,p):=\left(\int_{\mathbb N_0^n}|q-p|d\mu_n\right)^2,
	$$ 
	while the posterior predictive distribution $\big(P_{k,n}^*\big)^{P^n}$ is the Bayes estimator of the sampling distribution $P_\theta^n$ for the squared total variation loss function on $\mathbb N_0^n$. 
	
	Moreover, the image $M_n^*(k,\cdot):=\left[\big(P_{k,n}^*\big)^{P^n}\right]^{I_1}=I_1\big(P_{k,n}^*\big)^{P^n}$ is the Bayes estimator of the probability distribution $P_\theta$ for the squared total variation on $\mathbb N_0$ and its probability function
	$$k'\ge0\,\,\longmapsto\,\,\frac{dM_n^*(k,\cdot)}{d\mu_1}(k')=\int_0^\infty p_\theta(k')\cdot p^*_{k,n}(\theta)\,d\theta=\frac{(k'+\|k\|_1)!}{k'!\cdot (\|k\|_1)!}\cdot\frac{(\lambda+n)^{\|k\|_1+1}}{(\lambda+n+1)^{k'+\|k\|_1+1}}
	$$
	is the Bayes estimator of the probability function $p_\theta$ for the loss function $W'_1$.$\quad\Box$
\end{myexa}

\begin{myexa}\rm Let $P_\theta$ be the Bernoulli distribution with parameter $\theta\in(0,1)$ whose probability function is $p_\theta(k):=\theta^k(º-\theta)^{n-k}$, $k=0,1$. So $P_\theta^n$ is the joint distribution of a sample of size $n$ of a Bernoulli distribution with parameter $\theta$ and its probability function is
	$$p_{\theta,n}(k)=\theta^{\|k\|_1}(1-\theta)^{n-\|k\|_1},\quad k\in\{0,1\}^n
	$$
	where $\|k\|_1:=\sum_{i=1}^kk_i$. Consider the uniform distribution on the unit interval as prior distribution.
	So, the posterior distribution given  $k\in\{0,1\}^n$ is the Beta distribution 
	$$P^*_{k,n}=B(\|k\|_1+1,n-\|k\|_1+1)
	$$
	with parameters $\|k\|_1+1$ and $n-\|k\|_1+1$. Hence, denoting $\mu_n$ for the counter measure on $\{0,1\}^n$ and $\beta$ the Euler beta function,	
	the probability function of the posterior predictive probability given $k\in\{0,1\}^n$ is
	\begin{gather*}\begin{split}
			\frac{d(P_{k,n}^*)^{P^n}}{d\mu_n}(k')&=\int_\Theta p_{\theta,n}(k')\cdot p_{k,n}^*(\theta)\,d\theta\\&=\frac{\beta(\|k\|_1+\|k'\|_1+1)\beta(2n-\|k\|_1- \|k'\|_1+1)}{\beta(\|k\|_1+1)\beta(n-\|k\|_1+1)}\\
			&=\frac{\Gamma(n+2)}{\Gamma(2n+2)}\cdot
			\frac{(\|k'\|_1+\|k\|_1)!\cdot (2n-\|k'\|_1-\|k\|_1)!}{(\|k\|_1)!\cdot (n-\|k\|_1)!}.
	\end{split}\end{gather*}
	
	This is the Bayes estimator of the joint probability function $p_{\theta,n}$ for the loss function $W'_n(q,p):=\left(\int_{\{0,1\}^n}|q-p|d\mu_n\right)^2$, 	while the posterior predictive distribution $\big(P_{k,n}^*\big)^{P^n}$ is the Bayes estimator of the sampling distribution $P_\theta^n$ for the squared total variation loss function on $\{0,1\}^n$.
	
	Moreover, the image $M_n^*(k,\cdot):=\left[\big(P_{k,n}^*\big)^{P^n}\right]^{I_1}=I_1\big(P_{k,n}^*\big)^{P^n}$ is the Bayes estimator of the probability distribution $P_\theta$ for the squared total variation on $\{0,1\}$ and its probability function
	\begin{gather*}\begin{split}
	k'\in\{0,1\}\,\,\longmapsto\,\,&\frac{dM_n^*(k,\cdot)}{d\mu_1}(k')=\int_0^\infty p_\theta(k')\cdot p^*_{k,n}(\theta)\,d\theta\\&
	=\frac{\Gamma(n+2)}{\Gamma(2n+2)}\cdot
	\frac{(k'+\|k\|_1)!\cdot (2n-k'-\|k\|_1)!}{(\|k\|_1)!\cdot (n-\|k\|_1)!}
	\end{split}\end{gather*}
	is the Bayes estimator of the probability function $p_\theta$ for the $L^1$-squared loss function $W'_1$.$\quad\Box$
\end{myexa}

\section{Proofs}\par

\begin{myproo}{\sc (of Theorem 1)} \rm  (i)  
Notice that, writing $f_A(\theta):=P_\theta(A)$, 
$$\big(P_\omega^*\big)^P(A)=\int_\Theta P_\theta(A)dP_\omega^*(\theta)=E_{P_\omega^*}(f_A),
$$
that, as it is well known (see Nogales (1998) p. 147), is the Bayes estimator of $f_A$ for the squared error loss function. 

In the same way, if $X$	is a real integrable statistic on $(\Omega,\mathcal A)$ and $f(\theta):=E_\theta(X)$, we have that
$$E_{(P_\omega^*)^P}(X)=\int_\Theta\int_\Omega X(\omega')dP_\theta(\omega')dP_\omega^*(\theta)=E_{P_\omega^*}(f)
$$
is the Bayes estimator of $f$, the mean of $X$.

(ii) 	According (i), given $A\in\mathcal A$, 
$$\int_{\Omega\times\Theta} \left|\big(P_\omega^*\big)^P(A)-P_\theta(A)\right|^2d\Pi(\omega,\theta)\le \int_{\Omega\times\Theta} \left|X(\omega)-P_\theta(A)\right|^2d\Pi(\omega,\theta),
$$
for any real measurable function $X$ on $(\Omega,\mathcal A)$. If $\mathcal A$ is a separable $\sigma$-field, there exists a countable algebra  $\mathcal A_0$ such that $\mathcal A=\sigma(\mathcal A_0)$. In particular, it follows that
$$\sup_{A\in\mathcal A}\left|M(\omega,A)-P_\theta(A)\right|^2=
\sup_{A\in\mathcal A_0}\left|M(\omega,A)-P_\theta(A)\right|^2
$$
is $(\mathcal A\otimes\mathcal T)$-measurable. Given $(\omega,\theta)\in\Omega\times\Theta$, let 
$$C(\omega,\theta):=\sup_{A\in\mathcal A}\left|\big(P_\omega^*\big)^P(A)-P_\theta(A)\right|^2
$$
and, given $n\in\mathbb N$, choose $A_n\in\mathcal A_0$ so that 
$$C-\frac1n\le \left|\big(P_\omega^*\big)^P(A_n)-P_\theta(A_n)\right|^2.
$$
It follows from this  that 
\begin{gather*}\begin{split}
		\int_{\Omega\times\Theta} Cd\Pi&\le 
		\int_{\Omega\times\Theta} \left|\big(P_\omega^*\big)^P(A_n)-P_\theta(A_n)\right|^2d\Pi(\omega,\theta)+\frac1n\\
		&\le 		\int_{\Omega\times\Theta}\sup_{A\in\mathcal A}|M(\omega,A)-P_\theta(A)|^2d\Pi(\omega,\theta)+\frac1n,
\end{split}\end{gather*}
and this gives the proof as $n$ is arbitrary. To refine the proof from a measure-theoretical point of view, a judicious use of the Ryll-Nardzewski and Kuratowski measurable selection theorem would also be helpful. 
See the details in the remark at the end of the section. 

(iii) 	It follows from (ii) that, to estimate the density $p_\theta$, the posterior predictive density 
$$b^*_{Q,\omega}(\omega'):=\frac{d(P_\omega^*)^P}{d\mu}(\omega')$$
minimizes the Bayes mean risk for the loss function
$$W'_1(q,p):=\big(\textstyle\int|q-p|\,d\mu\big)^2,
$$
i.e.,
$$E_\Pi\left[\left(\int|b^*_{Q,\omega}-p_\theta|\,d\mu\right)^2\right]\le E_\Pi\left[\left(\int|m(\omega,\cdot)-p_\theta|\,d\mu\right)^2\right]
$$
for any measurable function $m:\Omega\times\Omega\rightarrow[0,\infty) $ such that $\int_\Omega m(\omega,\omega')d\mu(\omega')=1$ for every $\omega$. 		$\quad\Box$\end{myproo}

\begin{myproo} {\sc (of Theorem 2)} \rm (i)  
	Given $A\in\mathcal A^n$, Theorem 1.(i) shows that  the posterior predictive probability 
$\big(P_{\omega,n}^*\big)^{P^n}(A)$ of $A$ 
is the Bayes estimator of $f_A(\theta):=P_\theta^n(A)$ in the product Bayesian statistical experiment, as
$$\big(P_{\omega,n}^*\big)^{P^n}(A)=\int_\Theta P_\theta^n(A)dP_{\omega,n}^*(\theta)=E_{P_{\omega,n}^*}(f_A),
$$
i.e.
$$\int_{\Omega^n\times\Theta}\big|\big(P_{\omega,n}^*\big)^{P^n}(A)-P_\theta^n(A)\big|^2d\Pi_n(\omega,\theta)\le
\int_{\Omega^n\times\Theta}\big|X(\omega)-P_\theta^n(A)\big|^2d\Pi_n(\omega,\theta)
$$
for any other estimator $X:(\Omega^n,\mathcal A^n)\rightarrow\mathbb R$ of $f_A$. In particular, given $A\in\mathcal A$, applying this result to $I_1^{-1}(A)=A\times\Omega^{n-1}\in\mathcal A^n$, we obtain that
$$\int_{\Omega^n\times\Theta}\big|\big(P_{\omega,n}^*\big)^{P^n}(I_1^{-1}(A))-P_\theta(A)\big|^2d\Pi_n(\omega,\theta)\le
\int_{\Omega^n\times\Theta}\big|X(\omega)-P_\theta(A)\big|^2d\Pi_n(\omega,\theta)
$$
for any other estimator $X:(\Omega^n,\mathcal A^n)\rightarrow\mathbb R$ of $g_A:=P_\theta(A)$.

(ii) Being  $\mathcal A$ a separable $\sigma$-field, there exists a countable algebra  $\mathcal A_0$ such that $\mathcal A=\sigma(\mathcal A_0)$. In particular, it follows that
$$\sup_{A\in\mathcal A}\left|M(\omega,A)-P_\theta(A)\right|^2=
\sup_{A\in\mathcal A_0}\left|M(\omega,A)-P_\theta(A)\right|^2
$$
is $(\mathcal A\otimes\mathcal T)$-measurable. Given $(\omega,\theta)\in\Omega^n\times\Theta$, let 
$$C_n(\omega,\theta):=\sup_{A\in\mathcal A}\left|\big(P_{\omega,n}^*\big)^{P^n}(I_1^{-1}(A))-P_\theta(A)\right|^2
$$
and, given $k\in\mathbb N$, choose $A_k\in\mathcal A_0$ so that 
$$C_n-\frac1k\le \left|\big(P_{\omega,n}^*\big)^{P^n}(I_1^{-1}(A_k))-P_\theta(A_k)\right|^2.
$$
It follows that 
\begin{gather*}\begin{split}
\int_{\Omega\times\Theta} C_nd\Pi_n&\le 
\int_{\Omega\times\Theta} \left|\big(P_{\omega,n}^*\big)^{P^n}(I_1^{-1}(A_k))-P_\theta(A_k)\right|^2d\Pi_n(\omega,\theta)+\frac1k\\
&\le 		\int_{\Omega\times\Theta}\sup_{A\in\mathcal A}|M(\omega,A)-P_\theta(A)|^2d\Pi_n(\omega,\theta)+\frac1k,
\end{split}\end{gather*}
for any Markov kernel $M:(\Omega^n,\mathcal A^n)\pt (\Omega,\mathcal A)$ and, being $k$ arbitrary, this proves that 
$$M_n^*(\omega,A):=\big(P_{\omega,n}^*\big)^{P^n}(I_1^{-1}(A))
$$
is the Bayes estimator of $f(\theta):=P_\theta$ for the squared total variation loss function in the Bayesian statistical experiment
$$(\Omega^n,\mathcal A^n,\{P_\theta^n\colon \theta\in(\Theta,\mathcal T, Q)\})
$$ 
corresponding to a $n$-sized sample of the observed distribution. See the remark below. 

(iii) Note that, given $A\in\mathcal A$, Fubini's theorem yields 
$$\big(P_{\omega,n}^*\big)^{P^n}(I_1^{-1}(A))=\int_\Theta P_\theta(A)dP_{\omega,n}^*(\theta)=\int_A\int_\Theta p_\theta(\omega')\cdot p_{\omega,n}^*(\theta)dQ(\theta)d\mu(\omega'),
$$
where $p_{\omega,n}^*$ denotes the posterior density with respect to the prior distribution $Q$. Hence, for $\omega\in\Omega^n$, the $\mu$-density of $M_n^*(\omega,\cdot)$ is
$$\frac{dM_n^*(\omega,\cdot)}{d\mu}(\omega')=\int_\Theta p_\theta(\omega')\cdot p_{\omega,n}^*(\theta)dQ(\theta),
$$
and this is the Bayes estimator of the sampling density $p_\theta$ for the loss function $W'_1$. $\quad\Box$\end{myproo}

\begin{myrem}\rm (A precision on measure-theorethical technicalities in the proofs of the previous results) We detail the proof of Theorem 1.(ii), being that of Theorem 2.(ii) (and even that of the last remark of Section 3) similar. 
It follows from Theorem 1.(i) that, given $(\omega,\theta)\in\Omega\times\Theta$, and writing 
$$C(\omega,\theta):=\sup_{A\in\mathcal A}\left|\big(P_\omega^*\big)^P(A)-P_\theta(A)\right|^2,
$$
we have that, given $n\in\mathbb N$, there exists  $A_n(\omega,\theta)\in\mathcal A_0$ so that 
$$C(\omega,\theta)-\frac1n\le \left|\big(P_\omega^*\big)^P(A_n(\omega,\theta))-P_\theta(A_n(\omega,\theta))\right|^2.
$$
To continue the proof we will use the Ryll-Nardzewski and Kuratowski measurable selection theorem as appears in Bogachev (2007), p. 36. 

Keeping the notations of this book, we make $(T,\mathcal M)=(\Omega\times\Theta,\mathcal A\otimes \mathcal T)$ and $X=\mathcal A_0$ (the countable field generating $\mathcal A$). Given $n\in\mathbb N$, let us consider the map $S_n:\Omega\times\Theta\rightarrow \mathcal P(X)$ defined by
$$S_n(\omega,\theta)=\left\{A\in\mathcal A_0\colon C(\omega,\theta)-\frac1n\le\left|\big(P_\omega^*\big)^P(A)-P_\theta(A)\right|^2\right\}
$$
We have that $\emptyset\ne S_n(\omega,\theta)\subset X$ and $S_n(\omega,\theta)$ is closed for the discrete topology on $\mathcal A_0$. Moreover, given an open set $U\subset\mathcal A_0$,
$$\{(\omega,\theta)\colon S_n(\omega,\theta)\cap U\ne \emptyset\}\in\mathcal A\otimes\mathcal T
$$
because, given $A\in\mathcal A_0$, 
$$\{(\omega,\theta)\colon S_n(\omega,\theta)\ni A\}=
\left\{(\omega,\theta)\colon C(\omega,\theta)-\left|\big(P_\omega^*\big)^P(A)-P_\theta(A)\right|^2\le\frac1n\right\}\in\mathcal A\otimes\mathcal T.
$$
So, according to the measurable selection theorem cited above, there exists a measurable map $s_n:(\Omega\times\Theta,\mathcal A\otimes\mathcal T)\rightarrow\mathcal (A_0,\mathcal P(\mathcal A_0))$ such that $s_n(\omega,\theta)\in S_n(\omega,\theta)$ for every $(\omega,\theta)$ or, which is the same,
$$C(\omega,\theta)-\frac1n\le \left|\big(P_\omega^*\big)^P(s_n(\omega,\theta))-P_\theta(s_n(\omega,\theta))\right|^2.
$$
It follows that 
\begin{gather*}\begin{split}
		\int_{\Omega\times\Theta} C(\omega,\theta)d\Pi(\omega,\theta)&\le 
		\int_{\Omega\times\Theta} \left|\big(P_\omega^*\big)^P(s_n(\omega,\theta))-P_\theta(s_n(\omega,\theta))\right|^2d\Pi(\omega,\theta)+\frac1n\\
		&\le 		\int_{\Omega\times\Theta}\sup_{A\in\mathcal A}|M(\omega,A)-P_\theta(A)|^2d\Pi(\omega,\theta)+\frac1n,
\end{split}\end{gather*}
which gives the proof as $n$ is arbitrary. $\quad\Box$
\end{myrem}

\section{Acknowledgements}
This paper has been supported by the Junta de Extremadura (Spain) under the grant Gr18016.
\vspace{1ex}

\section* {References:}

\begin{itemize}

\item Barra, J.R. (1971) Notions Fondamentales de Statistique Mathématique, Dunod, Paris. 

\item Bogachev, V.I. (2007), Measure Theory, Vol. II, Springer, Berlin.

\item Boldstat, W.M. (2004) Introduction to Bayesian Statistics, Wiley, New Jersey. 

\item Ferguson, T.S. (1983) Bayesian density estimation by mixtures of normal distributions in "Recent advances in statistics", pages 287–302. Academic Press, New York.

\item Florens, J.P., Mouchart, M., and Rolin, J.M. (1990) Elements of Bayesian Statistics, Marcel Dekker, New York.

\item Geisser, S. (1993) Predictive Inference: An Introduction, Springer Science+ Business Media, Dordrecht.

\item Gelman, A., Carlin, J.B., Stern, H.S., Dunson, D.B., Vehtari, A., Rubin, D.B. (2014) Bayesian Data Analysis, 3rd ed., CRC Press.

\item Ghosal, S., Vaart, A.v.d. (2017) Fundamentals of Noparametric Bayesian Inference, Cambridge University Press, Cambridge UK. 

\item Ghosh J.K., Ramamoorthy, R.V. (2003) Bayesian
 Nonparametrics, Springer, New York. 
 
 \item Lijoi, A., Prünster, I. (2010) Models beyond the Dirichlet Process, in "Bayesian Nonparametrics", ed. by Hjort, N.L., Holmes, C., Müller, P., Walker, S.G., Cambridge Series in Statistical and Probabilistic Mathematics, Cambridge.

\item Lo, A.Y. (1984) On a class of Bayesian nonparametric estimates. I. Density estimates. Ann. Statist., 12(1):351–357.

\item Nogales, A.G. (1998) Estadística Matemática, Servicio de Publicaciones de la Universidad de Extremadura, Cáceres, Spain. 

\item Rubin, D.B. (1984) Bayesianly justifiable and relevant frequency calculations for the applied statisticians, The Annals of Statistics, 12(4), 1151-1172.

\end{itemize}

\end{document}